\documentclass[oneside]{amsart}

\usepackage[letterpaper,body={12.6cm,18.5cm}, mag=1000]{geometry}
\usepackage{amssymb}
\usepackage{amsthm}
\usepackage{amscd}
\numberwithin{equation}{section}

\theoremstyle{plain}

\theoremstyle{definition}

\newcommand{\dlabel}[1]{\ifmmode \text{\ttfamily \upshape [#1] } \else
{\ttfamily \upshape [#1] }\fi \label{#1}}

\newcommand{\Z}{\operatorname{Z} }

\newcommand{\Aut}{\operatorname{Aut} }

\begin{document}

\baselineskip 14pt

\title{On automorphisms of finite $p$-groups}
\author{Manoj K.~Yadav}

\begin{center}
\bf J. Group Theory, Vol. 10 (2007), 859-866
\end{center}

\address{School of Mathematics, Harish-Chandra Research Institute 
Chhatnag Road, Jhunsi, Allahabad - 211 019, INDIA}

\email{myadav@mri.ernet.in}

%\thanks{2000 MSC. 20D45; 20D15}
%\thanks{Key Words. Finite $p$-group, Central automorphism, Camina pair}
\thanks{Research supported by DST (SERC Division), the Govt. of INDIA}

%\date{\today}

\begin{abstract}
Let $G$ be a finite $p$-group such that $x\Z(G) \subseteq x^G$ for all $x \in
G - \Z(G)$, where $x^G$ denotes the conjugacy class of $x$ in $G$. Then $|G|$
divides $|\Aut(G)|$, where $\Aut(G)$ is the group of all automorphisms of $G$.
\end{abstract}

\maketitle

\end{document}